\documentclass[12pt] {article}
\usepackage{amssymb,amsxtra,amsmath,amsthm}
\def \ZZ{{\mathbb{Z}}}
\def \QQ{{\mathbb{Q}}}
\def \RR{{\mathbb{R}}}
\def \CC{{\mathbb{C}}}
\def \FF{{\mathbb{F}}}

\def \OOO{{\mathcal{O}}}
\def \TTT{{\mathcal{T}}}

\def \pp{{\mathfrak p}}

\def\det{\mathop{\mathrm{det}}}

\def\Ell{\mathop{\mathrm{Ell}}}
\def\mod{\mathop{\mathrm{mod}}}

\def\vert{\mathop{\mathrm{vert}}}
\def\edge{\mathop{\mathrm{edge}}}
\def\star{\mathop{\mbox{\Large$*$}}}
\def\Star{\mathop{\mbox{\LARGE$*$}}}

\def\diag{\mathop{\mathrm{diag}}}

\def\Gamma{\varGamma}
\def\Delta{\varDelta}
\def\Cl{\mathop{\mathrm{Cl}}}

\begin{document}

\begin{center}
{\Large {\bf Elliptic points of the Drinfeld modular groups}}\\
\bigskip
{\tiny {\rm BY}}\\
\bigskip
{\sc A. W. Mason}\\
\bigskip
{\small {\it Department of Mathematics, University of Glasgow\\
Glasgow G12 8QW, Scotland, U.K.\\
e-mail: awm@maths.gla.ac.uk}}\\
\bigskip
{\tiny {\rm AND}}\\
\bigskip
{\sc Andreas Schweizer\footnote{Most of this paper was written while
the second author was working at the Institute of Mathematics at
Academia Sinica in Taipei, supported by grant 99-2115-M-001-011-MY2
from the National Science Council (NSC) of Taiwan. During the final
stage the second author was supported by ASARC in South Korea}}\\
\bigskip
{\small {\it Department of Mathematics,\\
Korea Advanced Institute of Science and Technology (KAIST),\\
Daejeon 305-701\\
SOUTH KOREA\\
e-mail: schweizer@kaist.ac.kr}}
\end{center}
\begin{abstract}
Let $K$ be an algebraic function field with constant field $\FF_q$.
Fix a place $\infty$ of $K$ of degree $\delta$ and let $A$ be the
ring of elements of $K$ that are integral outside $\infty$. We give
an explicit description of the elliptic points for the action of the
Drinfeld modular group $G=GL_2(A)$ on the Drinfeld's upper half-plane
$\Omega$ and on the Drinfeld modular curve $G\!\setminus\!\Omega$.
It is known that under the {\it building map\;} elliptic points are
mapped onto vertices of the {\it Bruhat-Tits tree} of $G$. We show how such
vertices can be determined by a simple condition on their stabilizers. Finally
for the special case $\delta=1$ we obtain from this a surprising free
product decomposition for $PGL_2(A)$.
\end{abstract}

\subsection*{Introduction}

Let $K$ be an algebraic function field of one variable with constant
field $\FF_q$, the finite field of order $q$, and let $\infty$ be a
fixed place of $K$ of degree $\delta$. Let $K_{\infty}$ be the
completion of $K$ with respect to $\infty$ and let $C_{\infty}$ be
the $\infty$-completion of an algebraic closure of $K_{\infty}$.
The set $\Omega = C_{\infty} \backslash K_{\infty}$ is often
referred to as {\it Drinfeld's upper half-plane}. We denote the ring of
all those elements of $K$ which are integral outside $\infty$ by
$A$. (The simplest examples are $K=\FF_q(t)$ and $A=\FF_q[t]$.) The
group $G=GL_2(A)$ plays a fundamental role [D] in the theory of {\it
Drinfeld modular curves}. For this reason we will call $G$ a {\it
Drinfeld modular group}. Drinfeld [D] has extended the classical
theory of modular curves to the function field setting. Here $\QQ,
\RR, \CC$ are replaced by $K, K_{\infty}, C_{\infty}$, respectively.
The roles of the {\it classical upper half-plane}, $\mathbb{H}$,
(in $\CC$) and the {\it classical modular group}, $SL_2(\ZZ)$,
are assumed by $\Omega$ and $G$, respectively. The group $G$
acts as a set of linear fractional transformations on $\Omega$.
\\ \\
Let $S$ be a subgroup of $G$. We say that elements
$\omega_1,\;\omega_2 \in \Omega$ are $S${\it-equivalent}
if and only if $\omega_1=s(\omega_2)$, for some $s\in S$.
For each subgroup $S$ of $G$ and $\omega \in \Omega$, let
$S_{\omega}$ denote the {\it stabilizer} of $\omega$ in $S$.
\\ \\
\noindent {\bf Definition. }
The element $\omega\in\Omega$ is called an {\it elliptic}
element of $S$ if $S_{\omega}$ is non-trivial, i.e. it does {\it not}
consist entirely of scalar matrices. It is clear that $S$ acts on
its set of elliptic elements, $E(S)$. We put $\Ell(S)=S\backslash
E(S)$ and refer to its elements as the {\it elliptic points} of $S$.
\\ \\
Elliptic points are very important for a number of reasons.
One of the purposes of Drinfeld's theory is to provide an
analytical description for the so-called {\it Drinfeld modular
curve}, $G\backslash \Omega$ and hence $S\backslash \Omega$, for
every finite index subgroup $S$. Of particular importance in this
regard is, for example, the {\it genus} of such a curve whose
evaluation usually depends on the {\it Hurwitz formula} [G, p.87].
This relates the genera of $G \backslash \Omega$ and $S \backslash
\Omega$ and contains ramification factors which are in part
determined by elliptic points.
\\ \\
\noindent For $SL_2(\ZZ)$,
 it is a classical result that every element of $\mathbb{H}=\{z \in \CC: \mathrm{Im}z > 0\}$ which is fixed by a non-scalar matrix is $SL_2(\ZZ)$-equivalent to one of $i,\rho \in \mathbb{H}$, where
 $i^2=-1$ and $\rho^2+\rho+1=0$. Moreover every element of finite order in $SL_2(\ZZ)$ lies in
 the stabilizer of one of these ``elliptic" elements. It follows then that $SL_2(\ZZ)$ has precisely two ``elliptic points". As we shall see the situation for Drinfeld modular groups is much more complicated.
 \\ \\
\noindent Our first principal result provides a precise description
of an elliptic element.
\\ \\
\noindent {\bf Theorem A.} \it Fix any
$\varepsilon\in\FF_{q^2}\setminus\FF_q$. An element $\omega \in \Omega$
is an elliptic element of $G$ if and only if
$$\omega=\frac{\varepsilon+s}{t}$$
for some $s,t\in A\;(t\neq 0)$, for which
$$(\varepsilon^q+s)(\varepsilon+s)=tt', \ \hbox{\it with } \ t' \in A.$$
 \noindent \rm It follows that $G$ has elliptic elements if and only if $\delta$ is odd.
  Every elliptic element $\omega\in\Omega$ lies in $\FF_{q^2}K \backslash K$. We deduce from Theorem A that the stabilizer $G_{\omega}$ of every elliptic $\omega\in\Omega$ is isomorphic to
 $\cong\FF_{q^2}^*$. We are also able to deduce that $|\Ell(G)|=L_K(-1)$, where
$L_K(u)$ is the {\it L-polynomial} of $K$ [St, Section 5.1]. These deductions are already known [G, p.50]. However our approach is much simpler than that of Gekeler.
Moreover, we derive more precise information and interesting applications.
\\ \\
\noindent  The Galois automorphism of $\FF_{q^2}/\FF_q$ extends to that of
$\FF_{q^2}K/K$ and gives rise to a {\it conjugate map},
$\omega\mapsto\overline{\omega}$, on $E(G)$. Many of our results depend
on whether or not $\omega$ and $\overline{\omega}$ are $G$-conjugate.
Of particular interest in this context is the subset of $\Ell(G)$ consisting
of all those points corresponding to elliptic elements $\omega$ for which
$\omega$ and $\overline{\omega}$ are $G$-equivalent. We are able to identify
this subset with a certain group of involutions and for this reason we denote
it by $\Ell(G)_2$. It turns out, rather surprisingly perhaps, that
$|\Ell(G)_2|$, as with $|\Ell(G)|$,  does not depend on $A$, i.e. is
independent of the particular choice of $\infty$. We investigate how
$\Ell(G)$ and $\Ell(G)_2$ are related.
\\ \\
\noindent Associated with the group $GL_2(K_{\infty})$ is its {\it
Bruhat-Tits building} which in this case is a {\it tree}, $\TTT$.
See [Se, Chapter II, Section 1]. From this $G$ inherits an action
on $\TTT$. Most of our results involve the well-known {\it building map}
$$\lambda : \Omega \longrightarrow \TTT.$$
\noindent See [G, p.41], [GR, p.37]. Our next principal result
elaborates on the way elliptic elements are mapped into $\TTT$ under
the building map. It is known that, if $\omega \in E(G)$, then
$\lambda(\omega)=v$, for some $v \in \vert(\TTT)$, and that $G_{\omega}\leq G_v$. As usual $G_v$
denotes the stabilizer of the vertex $v$ of $\TTT$ in $G$. It is
known [Se, Proposition 2, p.76] that $G_v$ is always {\it finite}.
We prove the following.
\\ \\
\noindent {\bf Theorem B.} \it Suppose that $\delta$ is odd.\\ \\
\noindent (a) Let
$ v \in \vert(\TTT)$. Then
$$ v = \lambda(\omega),\; for\; some\; \omega \in E(G),\;\; if\;\; and\;\; only\;\; if\;\; q^2-1\; divides\; |G_v|.$$
\noindent (b) Suppose that $\omega\in E(\Omega)$ and $\lambda(\omega)=v$.
\begin {itemize}
\item[(i)] If $\omega,\;\overline{\omega}$ are $G$-equivalent, then $$G_v \cong GL_2(\FF_q).$$
\item[(ii)] Otherwise,
$$ G_v = G_{\omega} \cong \FF_{q^2}^*.$$
\end{itemize}

\noindent \rm Let $\widetilde{v}$ denote the image in $\vert(G \backslash \TTT)$ of a vertex $v$ of $\TTT$.
We put $\widetilde{K}= \FF_{q^2}K$.
\\ \\
\noindent {\bf Theorem C.} \it If $\delta$ is odd, there exist bijections between the following sets
\begin{itemize}
\item[(i)] vertices $\widetilde{v}$ of $G\!\setminus\!\TTT$ such that $q^2 -1$ divides $|G_v|$;
\item[(ii)] conjugacy classes (in $G$) of cyclic subgroups of $G$ of order $q^2 -1$;
\item[(iii)] the orbits of the $Gal(\widetilde{K}/K)$-action on $\Ell(G)$.
\end{itemize}
In particular, among the uncountably many points of $G\!\setminus\!\Omega$
lying over any given vertex $\widetilde{v}$ of $G\!\setminus\!\TTT$ there
are exactly
\begin{itemize}
\item one elliptic point if $G_v\cong GL_2(\FF_q)$;
\item two ($Gal(\widetilde{K}/K)$-conjugate) elliptic points if
$G_v\cong\FF_{q^2}^*$;
\item no elliptic points in all other cases.
\end{itemize}
\rm

\noindent \rm Finally we focus our attention on the important
special case where $\delta=1$. It can be shown that a vertex $v$ of
$\TTT$ gives rise to an {\it isolated} vertex of $G \backslash \TTT$
when (and {\it only} when) $\delta=1$ and $v=\lambda(\omega)$ as in
Theorem B. Isolated vertices are important for the following reason.
If such a vertex and its incident edge arise from a vertex $v$ and
incident edge $e$ of $\TTT$, then, from Bass-Serre theory [Se,
Theorem 13, p.55],
$$G \cong H \star_{\quad L }K,$$
\noindent where $H=G_v$ and $L=G_e$, the stabilizer of $e$.
\noindent Our final principal result is the following.
\\ \\
\noindent {\bf Theorem D.} \it Suppose that $\delta=1$. Then there
exists a subgroup $P$ such that
$$PGL_2(A)\cong\left(\Star_{i=1}^r\ZZ/(q+1)\ZZ\right)\star P.$$
Moreover, if $q\geq 8$ is fixed, then $r$ grows exponentially with
the genus of $K$.
\\ \\
\noindent \rm This decomposition has a number of interesting consequences.\\ \\
\noindent \rm We will use the following list throughout this paper.

\subsection*{Notation}

\begin{tabular}{ll}
$\FF_q$ & the finite field of order $q$;\\
$K$         & an algebraic function field of one variable with
constant field $\FF_q$;\\
$g(K)$         & the genus of $K$;\\
$L_K(u)$ & the $L$-polynomial of $K$;\\
$\infty$    & a chosen place of $K$;\\
$\delta$    & the degree of the place $\infty$;\\
$A$ & the ring of all elements of $K$ that are integral outside $\infty$;\\
$\widetilde{K}$  &  the quadratic constant field extension $\FF_{q^2}K$ of $K$;\\
$\widetilde{A}$  &  $\FF_{q^2}A$, the integral closure $A$ in $\widetilde{K}$;\\
$\nu$ & the additive, discrete valuation of $K$ defined by $\infty$;\\
$\pi$& a local parameter at $\infty$ in $K$;\\
$K_{\infty}$ & $\cong\FF_{q^{\delta}}((\pi))$, the completion of $K$ with
respect to $\infty$;\\
$\OOO_{\infty}$ & $\cong\FF_{q^{\delta}}[[\pi]]$, the valuation ring of $K_\infty$;\\
$C_\infty$ & the completion of an algebraic closure of $K_\infty$;\\
$\Omega$ & $=C_\infty -K_\infty$, Drinfeld's upper half-plane;\\
$\TTT$ & the Bruhat-Tits tree of $GL_2(K_\infty)$;\\
$G$ & the group $GL_2(A)$;\\
$G_w$ & the stabilizer in $G$ of $w \in \vert(\TTT) \cup \edge(\TTT)$;\\
$G_\omega$ & the stabilizer in $G$ of $\omega \in \Omega$;\\
$Z$ & the centre of $G$;\\
$\Cl(R)$ & the ideal class group of the Dedekind ring $R$;\\
$\Cl^0(F)$  & the divisor class group of degree $0$ of the function field $F$;\\
\end{tabular}
\\ \\ \\
\noindent We recall that $A$ is an arithmetic Dedekind domain with
$A^*=\FF_q^*$. In addition $\nu(a)\leq 0$, for all $a \in A$. Moreover
$\nu(a)=0$ if and only if $a \in \FF_q^*$. By definition $Z$
consists of all the scalar matrices $\alpha I_2$, where
$\alpha \in\FF_q^*$.
As usual, the {\it degree} of a prime ideal of $A$ or of a prime divisor of $K$
is the degree of its residue field over the constant field. By linear extension
one obtains the degree of any ideal or divisor.
\\ \\
It is well known that if $\delta$ is odd, the place $\infty$ of $K$ has exactly one
extension to $\widetilde{K}$, denoted by $\infty'$. In this case, $\widetilde{A}$
is the ring of all those elements of $\widetilde{K}$ which are integral outside
$\infty'$. We note that, if $\varepsilon \in \FF_{q^2}\backslash \FF_q$, then
$\widetilde{A}=A+\varepsilon A$. The action of $Gal(\widetilde{K}/K)$ on
$\widetilde{K}$ is given by
$$\overline{a+\varepsilon b}=a+\varepsilon^q b,$$
\noindent where $a,b \in K$. Also, for any set $J$ in $\widetilde{K}$, for example
if $J$ is an ideal of $\widetilde{A}$, we write $\overline{J}$ for the conjugate
set $\{\overline{x}\ :\ x\in J\}$.
\\

\subsection*{1. Elliptic elements on the Drinfeld upper halfplane $\Omega$}

\noindent \rm Before our first principal result we record some
elementary properties of non-trivial elements of elliptic point stabilizers.
\\ \\
{\bf Lemma 1.1.} \it Let $\omega \in \Omega$ be an
elliptic element and let
$M=\left[\begin{array}{cc}a&b\\c&d\end{array}\right]$ be a
non-scalar element of $G_{\omega}$. Then the minimal polynomial of
$\omega$ over $K$ is
$$ m_{\omega}(x)=x^2+\sigma x+\tau,$$
where $\sigma=(d-a)/c$ and $\tau=-b/c$.
\\ \\
{\bf Proof.}
\rm Follows from the fact that $M(\omega)=\omega$. Note that $bc\neq
0$, since $\omega \notin K$.
\hfill $\Box$
\\ \\ \noindent Before proceeding the following observation is critical.\\ \\
\noindent {\it The matrix $M \in G$ fixes $\omega \in \Omega$ if and only if $\left[\begin{array}{c}\omega\\1\end{array}\right]$ is an eigenvector of $M$.} \\ \\
{\bf Lemma 1.2.} \it
Let $\omega\in\Omega$ be an elliptic point and let $M\in G_{\omega}$
be non-scalar. Then $\omega\in\widetilde{K}$, and
$\left[\begin{array}{c}\omega\\1\end{array}\right]$
is an eigenvector of $M$ with eigenvalue
$\varepsilon\in\FF_{q^2}\setminus\FF_q$.
\\ \\
{\bf Proof.} \rm Let
$$M=\left[\begin{array}{cc}a&b\\c&d\end{array}\right].$$
Then $bc\neq 0$ by Lemma 1.1. It follows that $K(\omega)$ is a
quadratic extension of $K$. Now there exists $\varepsilon$ such that
$$ a\omega+b=\varepsilon\omega\ \ \mathrm{and}\ \ c\omega+d=\varepsilon.$$
Obviously $K(\omega)=K(\varepsilon)$. Moreover, $\varepsilon$ is an eigenvalue
of $M$ and so
$$\varepsilon^2+\eta\varepsilon+\rho=0,$$
where $\eta=-(a+d)$ and $\rho=\det(M)=(ad-bc) \in \FF_q^*$.
Let $B$ denote the integral closure of $A$ in $K(\varepsilon)$. Since
$M^{-1}\in G_{\omega}$ has eigenvalue $\varepsilon^{-1}$, we have
$\varepsilon,\varepsilon^{-1} \in B^*$. Now $\varepsilon \notin K_{\infty}$
(since $\omega \notin K_{\infty}$), so the place $\infty$ has only
one extension $\infty'$ to $K(\varepsilon)$, and $B$ consists of the elements
that are integral outside $\infty'$. Since $\varepsilon$ is invertible at all
places outside $\infty'$, by the product formula it must also be
invertible at $\infty'$ and hence a constant. So $\varepsilon$ is
algebraic over $\FF_q$ and since it generates a quadratic extension
of $K$ we conclude that $\varepsilon \in\FF_{q^2}\backslash \FF_q$. Thus
$$K(\omega)=K(\varepsilon)=\widetilde{K}.$$
\hfill $\Box$
\\ \\
We proceed to determine the stabilizer of an elliptic element.
\\ \\
{\bf Proposition 1.3.} \it Let $\omega$ be any elliptic element
of any $G$. Then $$G_{\omega} \cong \FF_{q^2}^*.$$
This isomorphism is given by mapping $M\in G_{\omega}$ to the eigenvalue
of $\left[\begin{array}{c}\omega\\1\end{array}\right]$ and it also
respects addition of matrices.
\\ \\
{\bf Proof.} \rm By Lemma 1.2
$\left[\begin{array}{c}\omega\\1\end{array}\right]$ is an eigenvector
for all $M \in G_{\omega}$ with corresponding eigenvalue
$\varepsilon \in \FF_{q^2}^*$ depending on $M$.
Applying $\tau\in Gal(\widetilde{K}/K)$, we see that
$\left[\begin{array}{c}\overline{\omega}\\1\end{array}\right]$
is an eigenvector with eigenvalue $\varepsilon^q$.
Hence there exists a matrix $X \in GL_2(\widetilde{K})$,
such that, for all $M \in G_{\omega}$,
$$ XMX^{-1}=\mathrm{diag}(\varepsilon,\varepsilon^q).$$
There is therefore a monomorphism
$$ G_{\omega} \hookrightarrow \FF_{q^2}^*.$$
\noindent To show that this map is surjective, we observe that by
definition $G_{\omega}$ contains a nonscalar $N$ with eigenvalues
$\mu,\mu^q \in\FF_{q^2}\setminus\FF_q$ and that for
all $\alpha,\beta \in \FF_q$, with $(\alpha,\beta) \neq (0,0)$,
$$ Y=\alpha I_2+\beta N \in G_{\omega},$$
and
$$XYX^{-1}=\mathrm{diag}(\alpha+\beta\mu,\alpha+\beta\mu^q).$$
\noindent The result follows.
\hfill $\Box$
\\ \\
For an alternative proof of Proposition 1.3 see [G, p.50].
\par
We are now able to provide a precise description of the elliptic
points of $G$.
\\ \\
{\bf Theorem 1.4.} \it Fix any $\varepsilon\in\FF_{q^2}\setminus\FF_q$. An
element $\omega \in \Omega$ is an elliptic element of $G$ if and only if
$$\omega=\frac{\varepsilon+s}{t}$$
for some $s,t\in A\;(t\neq 0)$, for which
$$(\varepsilon^q+s)(\varepsilon+s)=tt', \ \hbox{\it with } \ t' \in A.$$

\noindent {\bf Proof.} \rm Suppose $\omega \in \Omega$ is of the
form $\omega=\frac{\varepsilon +s}{t}$ as above. Let
$$M_0=\left[\begin{array}{cc}s'&-t'\\t&-s\end{array}\right],$$
where $s'=(\varepsilon+\varepsilon^q)+s$. Then it is easily verified that
(non-scalar) $M_0\in G_{\omega}$. Moreover $M_0$ has eigenvalues
$\varepsilon$ and $\varepsilon^q$ and determinant $\varepsilon^{q+1}\in\FF_q^*$.
\par
Conversely, let $\omega\in\Omega$ be elliptic. By Proposition 1.3 we
can choose $M=\left[\begin{array}{cc}a&b\\c&d\end{array}\right]$ in
$G_{\omega}$ such that $\varepsilon(M)=\varepsilon$. Then from the proof of
Lemma 1.2
$$\omega=(\varepsilon-d)/c.$$
\noindent Now $M(\omega)=\omega$ and so $$
c\omega^2+(d-a)\omega-b=0.$$ \noindent Let $\omega'$ be the other
root of this quadratic equation. Then $\omega'=\frac{\varepsilon^q
+s}{t}$ and
$$\omega\omega'=-b/c=(\varepsilon^q-d)(\varepsilon-d)/c^2.$$
Thus the condition is satisfied with $s=-d$ and $t=c$. \hfill $\Box$
\\ \\
{\bf Corollary 1.5.} \it $G$ has elliptic elements if
and only if $\delta$ is odd.
\\ \\
{\bf Proof.} \rm If $\omega$ is an elliptic element then, by definition,
$\omega \notin K_{\infty}$. By Theorem 1.4 there exists
$\varepsilon \in \FF_{q^2}\backslash \FF_q$ such that
$\varepsilon \notin K_{\infty}$. In addition,
$$\FF_{q^2} \subseteq K_{\infty}\; \Longleftrightarrow \;
\delta \;\mathrm{is\;even}.$$
On the other hand, if $\delta$ is odd, Theorem 1.4 implies that every
element of $\FF_{q^2} \backslash\FF_q$ is an elliptic point (s=0, t=1).
\hfill $\Box$
\\ \\
If $\omega=\frac{\varepsilon +s}{t}$ is elliptic and $M\in G_\omega$,
then from $M\omega=\omega$ one immediately obtains
$M\overline{\omega}=\overline{\omega}$.
So the conjugate $\overline{\omega}=\frac{\varepsilon^q +s}{t}$
is also elliptic with the same stabilizer, i.e.
$$G_{\overline{\omega}}=G_{\omega}.$$
A finer analysis of this in the next three sections will lead to
some interesting group-theoretic consequences. Among many others
we will need the following easy intermediate result.
\\ \\
{\bf Lemma 1.6.} \it If $\delta$ is odd, mapping $\{\omega, \overline{\omega}\}$
to $G_{\omega}=G_{\overline{\omega}}$ is a natural bijection between the unordered
pairs $\{\omega, \overline{\omega}\}$ of conjugate elliptic points and cyclic
subgroups of $G$ of order $q^2 -1$.
\\ \\
{\bf Proof.} \rm The inverse map is given by mapping the cyclic subgroup
to its two fixed points $\{\omega, \overline{\omega}\}$. These are indeed
elliptic. If not, they would lie in $K_\infty$, and consequently the eigenvalue
of the eigenvector $\left[\begin{array}{c}\omega\\1\end{array}\right]$
would be in $K_\infty\cap\FF_{q^2}^*=\FF_q^*$, in contradiction to the order
of the subgroup.
\hfill $\Box$
\\ \\
Lemma 1.6 also shows that if $\delta$ is odd, then the intersection of any
two cyclic subgroups of $G$ of order $q^2 -1$ is exactly $Z$.
\\ \\
\noindent We conclude this section with a further restriction on the
factor $t$ in Theorem 1.4 which we make use of later on.
\\ \\
{\bf Lemma 1.7.} \it Let $\omega=\frac{\varepsilon+s}{t}\in\Omega$ be an
elliptic element as in Theorem 1.4. Then
\begin{itemize}
\item[(a)] $\deg(\pp)$ is even for every prime ideal $\pp$ of $A$ that
divides $tA$.
\item[(b)] $\nu(t)$ is even.
\end{itemize}

\noindent{\bf Proof.} \rm  (a) Let $\pp$ be a prime
ideal of $A$ of odd degree. Then $\pp$ is inert in $\widetilde{A}$.
Let $\widetilde{\pp}$ be the prime ideal in $\widetilde{A}$ above
$\pp$. If $\pp$ divides $(t)$ in $A$, then $\widetilde{\pp}$ divides
$(\varepsilon+s)(\varepsilon^q +s)$ in $\widetilde{A}$. Since
$\widetilde{\pp}$ is a prime ideal, it must divide one of the two
factors. Applying the Frobenius automorphism of $\widetilde{K}/K$,
it also divides the other factor. Hence $\widetilde{\pp}$ divides
$(\varepsilon-\varepsilon^q)=\widetilde{A}$, a contradiction.

\noindent (b) By (a) and the product formula $\delta\nu(t)$ is even,
and $\delta$ is odd by Corollary 1.5.
\hfill $\Box$
\\

\subsection*{2. Elliptic points on the Drinfeld modular
curve $G\!\setminus\!\Omega$}

\noindent
{\it In view of Corollary 1.5 we assume throughout this section that
$\delta$ is odd.}\\

\noindent Central to the definition of $G \backslash \Omega$ is the
following equivalence relation.
\\ \\
\noindent {\bf Definition.} Let $\omega_1, \omega_2 \in \Omega$. We
say the $\omega_1,\omega_2$ are $G${\it -equivalent}, written
$\omega_1 \equiv \omega_2$, if and only if there exists $g \in G$
such that
$$\omega_1=g(\omega_2).$$
\noindent If $\omega_1=g(\omega_2)$ then
$$gG_{\omega_2}g^{-1}=G_{\omega_1}.$$
\noindent It follows that $G$-equivalent points of $\Omega$ have
isomorphic stabilizers in $G$. As we shall see the converse does not
hold. \noindent It is clear that $G$ acts on its elliptic points,
$E(G)$. We denote the set of equivalence classes by $\Ell(G)$. The
elements of this set are referred to as the {\it elliptic
points} of $G$.
\\ \\
\noindent In particular if $\delta$ is odd, then from Theorem 1.4
every $\varepsilon  \in \FF_{q^2}\setminus\FF_q$ is an elliptic point of
$G$. Moreover, if $\varepsilon$ and $\varepsilon'$ are any two elements of
$\FF_{q^2}\backslash \FF_q$, then $\varepsilon'=\alpha\varepsilon+\beta$ for
some $\alpha \in \FF_q^*,\; \beta \in \FF_q$ and hence
$$ \varepsilon \equiv \varepsilon'.$$
In particular, $\varepsilon \equiv \overline{\varepsilon}$. However, this
does not always hold for general elliptic points. For an arbitrary
elliptic point $\omega$ we will investigate later the precise
conditions under which $\omega$ and $\overline{\omega}$ are
$G$-equivalent. (They are not always equivalent despite the fact
that $G_{\omega}=G_{\overline{\omega}}$.)
\\ \\
\noindent {\bf Definition.} Let $A_0$ denote $A$ or $\widetilde{A}$.
If $I,I'$ are ideals in $A_0$, we write
$$ I \sim_{A_0} I' \Longleftrightarrow aI=bI',$$
for some non-zero $a,b \in A_0$. We use the standard notation
$\mathrm{Cl}(A_0)$ for the set of equivalence classes, usually
referred to as the {\it ideal class group of} $A_0$.
It is well-known that this is a {\it finite} group.
\\ \\
{\bf Lemma 2.1.} \it Let
$$\omega=\frac{\varepsilon+s}{t}$$
be an elliptic element,
where $\varepsilon,s,t$ are as defined in Theorem 1.4. Then
\begin{itemize}
\item[(a)] $J_{\omega}:=t A+(\varepsilon+s)A\; \trianglelefteq\; \widetilde{A}.$
\item[(b)] The ideal $J_{\omega}$ does not depend on the choice of
$\varepsilon\in\FF_{q^2}\setminus\FF_q$.
\end{itemize}

\noindent {\bf Proof.} \rm
(a) It suffices to prove that
$\varepsilon J_{\omega}\subseteq J_{\omega}$.
Now
$$ \varepsilon t=t(\varepsilon+s)-st \in J_{\omega}.$$
On the other hand,
$$\varepsilon(\varepsilon+s)=
(\varepsilon+\varepsilon^q)(\varepsilon+s)-(\varepsilon^q+s)
(\varepsilon+s)+s(\varepsilon+s) \in J_{\omega},$$
by the properties of $\varepsilon,s,t$.
\\
(b) Choosing a different $\varepsilon'\in\FF_{q^2}\setminus\FF_q$,
there exist $\alpha\in\FF_q^*$ and $\beta\in\FF_q$ with
$\varepsilon'=\alpha\varepsilon+\beta$. So
$\omega=\frac{\varepsilon'-\beta+\alpha s}{\alpha t}$,
which gives the same ideal.
\hfill $\Box$
\\ \\
\noindent Our next result is crucial since it enables us to identify
$\Ell(G)$ with a subgroup of $\mathrm{Cl}(\widetilde{A})$.
\\ \\
{\bf Lemma 2.2.} \it Let $\omega$ and $\omega'$ be elliptic elements
of $G$. Then
$$ \omega \equiv \omega' \;\Longleftrightarrow\;
J_{\omega}\sim_{\widetilde{A}}J_{\omega'}.$$

\noindent {\bf Proof.} \rm
Let $\omega =\frac{\varepsilon+s}{t}$ and
$\omega'=\frac{\varepsilon+s'}{t'}$. Then
$$tA+(\varepsilon+s)A\sim_{\widetilde{A}} t'A+(\varepsilon+s')A$$
if and only if there exist $a,b,c,d \in A$ with $ad-bc \in \FF_q^*$
and a non-zero $\rho\in\widetilde{K}$ such that
$$\rho t'=(at+b(\varepsilon+s))\;\
 \mathrm{and}\;\ \rho (\varepsilon+s')=(ct+d(\varepsilon+s)).$$
\hfill $\Box$
\\ \\
Mapping an elliptic element $\omega$ to the ideal class
$[J_{\omega}]\in\mathrm{Cl}(\widetilde{A})$ induces by Lemma 2.2
an injective map from $\Ell(G)$ into $\mathrm{Cl}(\widetilde{A})$.
In order to describe its image, we need the norm map $N$ from ideals
of $\widetilde{A}$ to ideals of $A$, and also from divisors of
$\widetilde{K}$ to divisors of $K$.
\par
If $\widetilde{P}$ is a prime ideal of $\widetilde{A}$,
then $N(\widetilde{P})=P^{f(\widetilde{P}/P)}$ where
$P=\widetilde{P}\cap A$ is the underlying prime ideal
of $A$ and $f(\widetilde{P}/P)$ is the inertia degree.
This definition is then canonically extended to products.
(See [ZS, Ch.~V, \S 11, p.306].)
Analogously for divisors (cf. [R2, pp.82]).
\par
In our simple situation we can equivalently say:
If $J\trianglelefteq\widetilde{A}$, then
$N(J)$ is the $A$-ideal $J\overline{J}\cap A$.
\par
Actually, $N(J)$ is also the $A$-ideal generated by all norms of
elements in $J$ [ZS, Ch.~V, \S 11, Lemma 3, p.307], but this
is not completely obvious. And in practice it is more awkward to
handle than the other properties.
\par
The norm map $N$ induces group homomorphisms
$$\overline{N}: \mathrm{Cl}(\widetilde{A})\longrightarrow\mathrm{Cl}(A)$$
and
$$\overline{N}: \mathrm{Cl}^0(\widetilde{K})\longrightarrow\mathrm{Cl}^0(K).$$
Our next goal is to show that the kernel of $\overline{N}$ is the
image of $\Ell(G)$. The following description of an element of
$\mathrm{Cl}(\widetilde{A})$ is essential for our purposes.
\\ \\
{\bf Lemma 2.3.} \it Let $J \trianglelefteq \widetilde{A}$.
Then, for any fixed $\varepsilon \in \FF_{q^2} \backslash \FF_q $,
there exist $a\in A$ and an ideal $I \trianglelefteq A$ such that
$$J \sim_{\widetilde{A}} J'=I+(\varepsilon+a) A.$$
\noindent Moreover $J' \cap A= I$ and $N(J')=I$.
\\ \\
{\bf Proof.} \rm Now $\widetilde{A}=A+\varepsilon A$ and so,
by [B, Chapter VII, Section 4.10, Proposition 24], there exists
$a,b \in A$ and an $A$-module $I'$, $A$-isomorphic to an $A$-ideal
such that
$$ J=I'+(a+\varepsilon b) A.$$
\noindent Since $A$ is Dedekind there are two possibilities.
\\ \\
\noindent \bf {(a)} $I'=x A,\; for\; some\; nonzero\;
x\in\widetilde{A}$:
\\ \\
\rm Then
$$ J \sim_{\widetilde{A}} \overline{x}J.$$
\noindent By multiplying by another term in $A$ (to ``clear denominators")
we may assume that $x \in A$.
\\ \\
\noindent {\bf (b)} $I'= Ax+Ay,\;  with\;
ey=fx\neq 0,\;where\;x,y \in \widetilde{A}\;and\; e,f \in A$:
\\ \\
\rm Replacing $J$ with $fy^{-1}J$ and then ``clearing denominators" as
above we may assume that $x,y \in A$.
\\ \\
\noindent From now on we replace $I'$ with $I$, where $I\trianglelefteq A$.
Let $i\in I$. Then $i\varepsilon \in J$ and so $i=bb'$, where $b' \in A$.
On the other hand $\varepsilon(a+\varepsilon b) \in J$, and since
$\varepsilon^2 =\alpha\varepsilon+\beta$ with $\alpha,\beta\in\FF_q$,
this implies $a=bb''$, where $b'' \in A$.
We now replace $J$ with $J'=b^{-1}J$, which has the desired form.
\par
Moreover, $J'\cap A= I$ is obvious. Finally,
$$J'\overline{J'}=I^2+I(\varepsilon-a)+I(\varepsilon^q -a)+
(\varepsilon-a)(\varepsilon^q -a)A
\subseteq I^2+I\widetilde{A}+(J'\cap A)\subseteq I\widetilde{A}.$$
Conversely, $J'\overline{J'}$ contains $I(\varepsilon-a)-I(\varepsilon^q -a)$,
and hence $I(\varepsilon-\varepsilon^q)=I\widetilde{A}$.
So together $J'\overline{J'}=I\widetilde{A}$ and thus $N(J')=I$.
\hfill $\Box$
\\ \\
\noindent {\bf Theorem 2.4.} \it
Mapping an elliptic element $\omega$ to the ideal class
$[J_{\omega}]$ in $\mathrm{Cl}(\widetilde{A})$ induces a bijection
between $\Ell(G)$ and the kernel of the surjective norm map
$\overline{N}: \mathrm{Cl}(\widetilde{A}) \longrightarrow \mathrm{Cl}(A)$.
\\ \\
\noindent {\bf Proof.} \rm
If $\omega=\frac{\varepsilon+s}{t}$ is elliptic, then the ideal
$J_{\omega}=tA+(\varepsilon+s)A$ from Lemma 2.1 has norm $tA$ by
Lemma 2.3. So $[J_{\omega}]$ lies in the kernel of $\overline{N}$.
\par
Conversely, we represent each element $[J]$ of $\mathrm{Cl}(\widetilde{A})$
by an ideal $J$ of the form given by Lemma 2.3. Then $[J] \in
\mathrm{Ker}\;\overline{N}$ if and only if $N(J)=I$ is principal,
i.e. if and only if
$$J=A c+A(a+\varepsilon),$$
for some non-zero $a,c \in A$.
\noindent Note that, if $J$ is of this form, then
$(a+\varepsilon)(a+\varepsilon^q)=cc'$, for some $c' \in A$, since
$J\cap A=I$. Suppose that $A c+A a\neq A$. Then there
exists a prime $\widetilde{A}$-ideal, $\pp$, containing $a,c$.
Thus $(a+\varepsilon)(a+\varepsilon^q) \in \pp$ and so $\varepsilon \in \pp$,
which implies that $\pp=\widetilde{A}$. Hence $A c+A
a=A$ and so $J$ is determined by the elliptic point
$\omega=(a+\varepsilon)/c$. (See Theorem 1.4.)
\par
Finally, we prove the surjectivity of
$\overline{N}: \mathrm{Cl}(\widetilde{A}) \longrightarrow \mathrm{Cl}(A)$.
Since $\delta$ is odd and hence $\infty$ is inert in $\widetilde{K}$
by [R2, Proposition 8.13], we can apply [R1, Proposition 2.2] which
tells us that $\overline{N}$ is surjective.
In passing we point out that we do not know whether there exists an
elementary proof of the surjectivity of $\overline{N}$, that is, a proof
that avoids the use of class field theory.
\hfill $\Box$
\\ \\
So far, $|\Ell(G)|$ seems to depend on the ring $A$, of which there are
infinitely many non-isomorphic ones in the same function field $K$.
But one can go one step further.
\\ \\
\noindent{\bf Lemma 2.5.} \it
The canonical map from $\mathrm{Cl}^0(\widetilde{K})$ to
$\mathrm{Cl}(\widetilde{A})$ restricts to an isomorphism of abelian
groups between the kernel of the surjective norm map
$\overline{N}: \mathrm{Cl}^0(\widetilde{K}) \longrightarrow \mathrm{Cl}^0(K)$
and the kernel of
$\overline{N}: \mathrm{Cl}(\widetilde{A}) \longrightarrow \mathrm{Cl}(A)$.
\\ \\
\noindent{\bf Proof.} \rm
Mapping the divisor $\prod\limits P^{e_P}$ of $\widetilde{K}$ to the
fractional ideal $\prod\limits_{P\neq\infty}P^{e_P}$ of $\widetilde{A}$
induces an isomorphism from $\Cl^0(\widetilde{K})$ to a subgroup of index
$\delta$ in $\Cl(\widetilde{A})$, namely to the classes consisting of
ideals whose degrees are divisible by $\delta$.
(Compare [R2, Proposition 14.1].)
But the degree of every principal ideal of $A$ obviously is divisible
by $\delta$. So if the ideal class $[J]$ is in the kernel of
$\overline{N}$, then $\delta$ divides $\deg(N(J))=2\deg(J)$ and hence
$\deg(J)$ since $\delta$ is odd.
Now one easily verifies that the map induces the desired isomorphism.
\par
As explained before, or by [R1, Lemma 1.2], we have
$$ |\mathrm{Cl}(\widetilde{A})|=\delta|\mathrm{Cl}^0(\widetilde{K})|\ \;
\mathrm{and}\ \;|\mathrm{Cl}(A)|=\delta|\mathrm{Cl}^0(K)|.$$
So the surjectivity of the norm map from $\mathrm{Cl}(\widetilde{A})$
to $\mathrm{Cl}(A)$ implies the surjectivity of the norm map from
$\mathrm{Cl^0}(\widetilde{K})$ to $\mathrm{Cl^0}(K)$.
\hfill $\Box$
\\ \\
{\bf Corollary 2.6.} \it With the above notation,
$$|\Ell(G)|=L_K(-1).$$

\noindent {\bf Proof.} \rm Combining Theorem 2.4 and Lemma 2.5 with
[St, Theorem V.1.15 (c),(f)], we have
$$|\Ell(G)|=\frac{|\mathrm{Cl^0}(\widetilde{K})|}{|\mathrm{Cl^0}(K)|}
=\frac{L_{\widetilde{K}}(1)}{L_K(1)}=L_K(-1).$$
\hfill $\Box$
\\ \\
\noindent Corollary 2.6 (as well as Proposition 1.3) is already known [G, p.50].
However our approach is much simpler than that of Gekeler. In particular it avoids
any mention of the fact that $G\backslash\Omega$ is a component of the moduli scheme
for Drinfeld $A$-modules of rank $2$. In addition, at this stage we don't yet need
the {\it building map} $\lambda : \Omega \longrightarrow
\TTT$, where $\TTT$ is the {\it Bruhat-Tits tree} associated with
$G$. (See [Se, Chapter II, Section 1.1], [G, p.41].)
The remaining results in this section will elaborate on the structure of $\Ell(G)$.
\\ \\
\noindent {\bf Lemma 2.7.} \it
\begin{itemize}
\item[(a)] For $q\geq 4$ there exists only one non-rational function
field $K$ with $|\Ell(G)|=1$, namely
$$K=\FF_4(x,y)\ \ \hbox{\it with}\ \ y^2 +y=x^3.$$
\item[(b)] More generally, for any positive integer $n$ there are only
finitely many nonrational function fields $K$ with $q\geq 3$ and
$|\Ell(G)|=n$.
\end{itemize}

\noindent {\bf Proof.} \rm
Using Corollary 2.6 and the Riemann Hypothesis for function fields
[St, Theorem 5.2.1], [St, Theorem 5.1.15(e)] we have
$$n=|\Ell(G)|=L_K(-1) \geq (\sqrt{q}-1)^{2g}.$$
For given $n$ this bounds $q$, and for $q>4$ it also bounds $g$.
\par
In particular, $n=1$ is only possible for $q\leq 4$; and if $n=1$
for $q=4$, then necessarily $L_K(u)=(1+2u)^{2g}$.
\par
A function field $K$ over $\FF_q$ with $L_K(u)=(1+\sqrt{q}u)^{2g}$
is called {\it maximal}. Equivalently, a maximal function field
is a function field with $q+1+2g\sqrt{q}$ places of degree $1$.
\par
By Ihara's Theorem [St, Proposition 5.3.3] the genus of a maximal
function field is bounded by $g\leq\frac{q-\sqrt{q}}{2}$.
For $q=4$ this leaves only the possibility $g=1$.
But it is well known that $y^2+y=x^3$ is the only elliptic function
field over $\FF_4$ with $L$-polynomial $(1+2u)^2$.
Alternatively one could invoke [RSt, Theorem] here. This finishes
the proof of (a).
\par
For (b) we still have to take care of the cases $q=3$ and $4$.
We exploit the following lower bound for the class number from
[St, Exercise 5.8, p.213]
$$L_K(1)\geq
\frac{q-1}{2}\cdot\frac{q^{2g}+1-2gq^g}{g(q^{g+1}-1)}\geq
\frac{q-1}{2}\cdot\frac{q^{2g}-2gq^g}{gq^{g+1}}=
\frac{q-1}{2}\left(\frac{q^{g-1}}{g}-\frac{2}{q}\right).$$
Applied to the field $\widetilde{K}$ this yields
$$L_{\widetilde{K}}(1)\geq \frac{c\cdot q^{2g}}{g}$$
where $c$ is a nonzero constant depending on $q$.
Combined with the upper bound
$$L_K(1)\leq (\sqrt{q}+1)^{2g}$$
from the Riemann Hypothesis [St, Theorem 5.2.1],
[St, Theorem 5.1.15(e)] this shows
$$|\Ell(G)|=\frac{L_{\widetilde{K}}(1)}{L_K(1)}\geq
\frac{c\cdot q^{2g}}{g(\sqrt{q}+1)^{2g}}.$$
So $|\Ell(G)|$ goes to infinity with $g$ provided $q\geq 3$.
\hfill $\Box$
\\ \\
We apply Lemma 2.2 to the cases for which $L_K(-1)=1$.
(One such is the genus zero case $K=\FF_q(T)$.)
Let $\varepsilon \in \FF_{q^2} \backslash \FF_q$.
Then, if $\omega$ is any elliptic point, there exists
$g \in G$ such that $g(\omega)=\varepsilon$.
In particular, then $\omega\equiv\overline{\omega}$
for all elliptic points.
\par
We will determine now when this happens in general.
\\ \\
By Lemma 2.3 we have
$$J_{\omega}J_{\overline{\omega}}=J_{\omega}\overline{J_{\omega}}=
N(J_{\omega})\widetilde{A}=t\widetilde{A}.$$
It follows that in $\mathrm{Cl}(\widetilde{A})$
$$[J_{\overline{\omega}}]=[J_{\omega}]^{-1}.$$
Hence
$$ \omega \equiv \overline{\omega}\Longleftrightarrow
[J_{\omega}]^2=1$$
in $\mathrm{Cl}(\widetilde{A})$, or equivalently, in
$\mathrm{Cl}^0(\widetilde{K})$.
\\ \\
\noindent {\bf Definition.}
Let $\Ell(G)_2$ be the subset of $\Ell(G)$ consisting of those orbits
of elliptic elements for which $\overline{\omega} \equiv \omega$.
\\ \\
We have just proved the following result.
\\ \\
\noindent {\bf Theorem 2.8.} \it
The bijection between $\Ell(G)$ and the kernel of the norm map
$\overline{N}$ described in Theorem 2.4 restricts to a bijection
between $\Ell(G)_2$ and the $2$-torsion subgroup of the kernel
of $\overline{N}$ in $\mathrm{Cl}(\widetilde{A})$, or by
Lemma 2.5 equivalently, the $2$-torsion subgroup of the kernel
of $\overline{N}$ in $\mathrm{Cl}^0(\widetilde{K})$.
\par
In particular, $|\Ell(G)|$ and $|\Ell(G)_2|$ only depend on $K$, not
on the choice of the place $\infty$ (apart from the general condition
that $\delta$ has to be odd).
\rm
\\ \\
Hence if $\Ell(G)=\Ell(G)_2$ (for example, when $L_K(-1)=1$)
it follows that
$\omega \equiv \overline{\omega}$
\noindent for all $\omega \in E(G)$. On the other hand we can prove
the following.
\\ \\
\noindent
{\bf Theorem 2.9.} \it
\begin{itemize}
\item[(a)] For $q\geq 8$ there are only two function fields $K$ of genus
$g>0$ for which $\Ell(G)=\Ell(G)_2$, namely
$$K=\FF_9(x,y)\ \ \hbox{\it with}\ \ y^3+y=x^4\ \ \hbox{\it (genus $3$)}$$
and
$$K=\FF_9(x,y)\ \ \hbox{\it with}\ \ y^2=x^3-x\ \ \hbox{\it (genus $1$)}.$$
\item[(b)] For fixed $q\geq 8$ we have
$\lim\limits_{g\to\infty}\frac{|\Ell(G)_2|}{|\Ell(G)|}=0$.
\end{itemize}

\noindent
{\bf Proof.} \rm (a) By Corollary 2.6 and the Riemann Hypothesis
for function fields [St, Theorem 5.2.1], [St, Theorem 5.1.15(e)]
$$|\Ell(G)|=L_K(-1) \geq (\sqrt{q}-1)^{2g}.$$
\noindent On the other hand, the $2$-torsion rank of an abelian
variety of dimension $g$ is bounded by $2g$, and even by $g$ if
the characteristic is $2$. Applying this to
$\mathrm{Cl}^0(\widetilde{K})$ (compare [R2, Chapter 11])
we get
$$|\Ell(G)_2| \leq 2^{2g},$$
and even $|\Ell(G)_2| \leq 2^g$ if the characteristic is $2$.
This proves both claims if $q>9$ and also if $q=8$.
\par
For the remaining case $q=9$ we note that by the same argument
$|\Ell(G)|=|\Ell(G)_2|$ is only possible if $L_K(u)=(1+3u)^{2g}$,
that is, if $K$ is a maximal function field. Then Ihara's Theorem
[St, Proposition 5.3.3] implies $g\leq 3$. Moreover, $g=2$ is not
possible, because then $K$ would be hyperelliptic, i.e. a double
covering of a rational function field $\FF_9(T)$, and hence could
have at most $2(9+1)<22$ places of degree $1$.
\par
By [RSt, Theorem] there is a unique maximal function field of
genus $3$ over $\FF_9$, namely the {\it Hermitian function field}
$\FF_9(x,y)$ with $y^3 +y=x^4$. Furthermore, by [RSt, Lemma 1] this
function field has $\mathrm{Cl}^0(K)\cong\bigoplus_{i=1}^{6}\ZZ/4\ZZ$.
Since $L_{\widetilde{K}}(t)=(1-9t)^6$, by the same argument we have
$\mathrm{Cl}^0(\widetilde{K})\cong\bigoplus_{i=1}^{6}\ZZ/8\ZZ$, and
hence the kernel of the norm map is indeed isomorphic to
$\bigoplus_{i=1}^{6}\ZZ/2\ZZ$.
\par
For $g=1$ we use the well-known fact that $y^2=x^3-x$ is the only
elliptic function field over $\FF_9$ with $L$-polynomial $(1+3u)^2$
or some explicit calculations with Weierstrass equations.
\par
Finally, to prove claim (b) for $q=9$ we bound $|\Ell(G)|$ from
below by exactly the same procedure as in the proof of Lemma 2.7.
Then $\frac{|\Ell(G)_2|}{|\Ell(G)|}\leq
\frac{g\cdot 8^{2g}}{c\cdot 9^{2g}}$, which goes to $0$.
\hfill $\Box$
\\ \\
When $q>9$ and $g>0$ therefore $G$ has an elliptic point $\omega_0$
which is {\it not} equivalent to $\overline{\omega_0}$. As we shall
see in the next two sections, points like these have a special significance
for the Bruhat-Tits tree and the structure of $G$.
\\ \\
\noindent {\bf Remark 2.10.}
It is not clear whether for $q\leq 7$ there are only finitely many
function fields $K$ with $\Ell(G)=\Ell(G)_2$.
And even if one could prove finiteness, the actual determination
of all such fields would be a tedious task.
\par
Let us consider the special case where $K$ is a quadratic extension
of a rational function field $\FF_q(T)$, that is, $K$ is hyperelliptic
or possibly elliptic. In this case the degree $4$ Galois extension
$\widetilde{K}/\FF_q(T)$ has $3$ intermediate extensions, namely
$K$, $\FF_{q^2}(T)$, and the unramified quadratic twist of $K$, which
we denote by $K'$.
\par
Now the kernel of the norm map from
$\mathrm{Cl}^0(\widetilde{K})$ to $\mathrm{Cl}^0(K)$ is isomorphic
to $\mathrm{Cl}^0(K')$. So the determination of all hyperelliptic
$K$ with $\Ell(G)=\Ell(G)_2$ is equivalent to the determination of
all hyperelliptic function fields with divisor class group of
exponent $2$.
\par
For the more special case where in addition a degree $1$ place of
$\FF_q(T)$ is ramified in $K'$ this is the goal of the paper [BD].
But even then case-by-case arguments and a computer search were needed.
\par
More importantly, on the way from [BD, Theorem 21] to [BD, Theorem 37]
several cases, including among others for $h=8$ the cases $q=5$, $g=2$
and $q=3$, $g=3,4$ as well as $q=2$, $4\leq g\leq 8$ seem to have got
lost, and consequently the main result of that paper is incomplete.
Without claim for completeness we point out some missing elliptic
function fields $K'$ with
$\Cl^0(K')\cong\ZZ/2\ZZ\oplus\ZZ/2\ZZ$, to wit
\par
\begin{tabular}{ll}
\ \ \ $K'=\FF_3(x,y)$ \ with & $y^2=x^3-x$,\\
\ \ \ $K'=\FF_5(x,y)$ \ with & $y^2=x^3+x$,\\
\ \ \ $K'=\FF_7(x,y)$ \ with & $y^2=x^3-1$,\\
\ \ \ $K'=\FF_9(x,y)$ \ with & $y^2=x^3-\sqrt{-1}x$.\\
\end{tabular}
\\
The last example is the unramified quadratic twist of the exceptional
$K$ in Theorem 2.9 (a).
\\

\subsection*{3. The images of elliptic points on the Bruhat-Tits tree $\TTT$}

\noindent Associated with the group $GL_2(K_{\infty})$ is its
{\it Bruhat-Tits building} which in this case is a
$(q^{\delta}+1)$-regular {\it tree}, $\TTT$.
The most convenient description for our purposes is the one
in [Se, Chapter II, Section 1]. See also [GR, Section 1.3].
The vertices of $\TTT$ are the homothety classes of
$\OOO_\infty$-lattices of rank $2$ in $K_\infty\oplus K_\infty$.
Two such vertices are joined by an edge if they contain lattices
$L_1$ and $L_2$ such that $L_2$ is a maximal $\OOO_\infty$-sublattice
of $L_1$. This definition is of course symmetric, because then
$\pi L_1$ is a maximal sublattice of $L_2$.
\par
Via its natural embedding into $GL_2(K_\infty)$, the group $G$ acts
on $\TTT$ {\it without inversion} [Se, Corollary, p.75]. Classical
Bass-Serre theory [Se, Theorem 13, p.55] shows how the structure of
$G$ can be derived from that of the quotient graph $G\backslash
\TTT$. The structure of this quotient is described in [Se, Theorem
9, p.106]. (Serre's approach uses the theory of vector bundles. For
a more elementary approach see [M, Theorem 4.7].)
In the sequel we will write $v$ and $e$ for vertices respectively edges
of $\TTT$ and $\widetilde{v}$ and $\widetilde{e}$ for their images in
$G\backslash\TTT$.
\par
A central object in the study of Drinfeld's half-plane is the
{\it building map}
$$\lambda : \Omega \longrightarrow \TTT.$$
\noindent See [G, p.41], [GR, Section 1.5].
We only mention the facts that we need and refer to the literature
for a thorough description.
\par
If $|\cdot|$ denotes the multiplicative valuation on $C_\infty$, then
every $\omega\in\Omega$ defines a norm $v_\omega(u,v):=|u\omega +v|$
on the vector space $K_\infty \oplus K_\infty$. By a theorem of Goldman
and Iwahori there are two types of such norms. If the unit ball of
$v_\omega$ is an $\OOO_\infty$-lattice $L$ in $K_\infty \oplus K_\infty$,
then $\lambda(\omega)$ is the vertex of $\TTT$ given by  the homothety
class of $L$. In all other cases $v_\omega$ is a ``convex combination"
of two norms that are of the former type and belong to two neighbouring
vertices. Correspondingly $\lambda$ then maps $\omega$ to a point on the
edge joining these two vertices.
\par
Another important feature is that $\lambda$ respects the actions of
$GL_2(K_\infty)$ on $\Omega$ and $\TTT$, that is
$$\lambda(g(\omega))=g(\lambda(\omega)).$$
In particular, $\lambda$ induces a map from the quotient space
$G\!\setminus\!\Omega$ to the quotient graph $G\!\setminus\!\TTT$.
Important information about $G\!\setminus\!\Omega$ is encoded in
the (in a certain sense) simpler object $G\!\setminus\!\TTT$ (see
for example [GR]). Here we explore this theme with respect to
elliptic points.
\\ \\
{\bf Lemma 3.1.} \it
If $\omega\in\Omega$ is an elliptic element, then $\lambda(\omega)$ is a vertex
of $\TTT$ and $G_{\omega}$ is a subgroup of $G_{\lambda(\omega)}$.
Moreover, $\lambda(\omega)=\lambda(\overline{\omega})$.
\rm
\\ \\
{\bf Proof.}
If $\delta$ is odd, for every $\varepsilon\in\FF_{q^2}\setminus\FF_q$ the
associated norm $v_{\varepsilon}((u,v))=|u\varepsilon +v|$ on
$K_\infty \oplus K_\infty$ obviously is the maximum norm
$\max\{|u|,|v|\}$, whose unit ball is the standard lattice
$\OOO_\infty \oplus \OOO_\infty$. So all $\varepsilon\in\FF_{q^2}\setminus\FF_q$
map to the standard vertex in $\TTT$.
\par
Now if $\omega$ is any elliptic point, by Theorem 1.4 we have
$\omega=\frac{\varepsilon +s}{t}$ and
$\overline{\omega}=\frac{\varepsilon^q +s}{t}$ for suitable $s,t\in A$.
So under $\lambda$ both, $\omega$ and $\overline{\omega}$ map to the
same vertex of $\TTT$, namely the image of the standard vertex under the
action of ${1\ s\choose 0\ t}\in GL_2(K_\infty)$
\par
The fact $G_\omega\leq G_{\lambda(\omega)}$ is clear.
(See also [GR, (1.5.3), p.37].)
\hfill $\Box$
\\ \\
We recall [Se, Proposition 2, p.76] that the elements of
finite order in $G$ are precisely those in
$$\bigcup_{v \in\vert(\TTT)}G_v.$$
\noindent We note that
$$Z\leq G_e\cap G_\omega,$$
\noindent for all $e \in \edge(\TTT)$ and $\omega \in E(G)$. Hence
$q-1$ divides all $|G_v|$.
\\ \\
\noindent Let $\omega \in E(G)$. Then we know that $G_{\omega} \leq
G_{\lambda(\omega)}$ and consequently $q^2-1$ divides
$|G_{\lambda(\omega)}|$, by Proposition 1.3. One of the main aims of this section is to establish the converse of this result.
However, for that we need a few lemmata.
\\ \\
{\bf Lemma 3.2.} \it
Let $M\in G_v$. Then the eigenvalues of $M$ lie in $\FF_{q^2}$.
\rm
\\ \\
{\bf Proof.}
The characteristic polynomial of $M$ is
$$t^2-\tau t+\eta,$$
where $\tau=\mathrm{tr}(M)$ and $\eta=\det(M)\in\FF_q^*$.
Now $M$ has finite order and so $\tau$ lies in the algebraic closure
of $\FF_q$ in $A$ which is $\FF_q$.
\hfill $\Box$
\\ \\
\noindent \rm Our next result, although little more than an
observation, is crucially important.
\\ \\
{\bf Lemma 3.3.} \it Let $w \in \vert(\TTT)\cup\edge(\TTT)$.
Suppose that $M_1,M_2$ are matrices in $G_w$. If
$$\det(\alpha_1M_1+\alpha_2M_2) \in \FF_q^*,$$
where $\alpha_1,\alpha_2 \in \FF_q$, then
$$\alpha_1M_1+\alpha_2M_2 \in G_w.$$

\noindent {\bf Proof.} \rm If $M_i$ fixes a vertex, that is,
a lattice class $\Lambda$, then because of $M_i \in GL_2(A)$
by [Se, II.1.3 Lemma 1, p.76] it fixes any underlying lattice
$L$. Thus $\alpha_1 M_1 +\alpha_2 M_2$ is an endomorphism of
$L$. But since its determinant is invertible in $\OOO_{\infty}$,
it actually is an automorphism of $L$. So it fixes the same
lattice class.
\hfill $\Box$
\\ \\
\noindent As is clear from the proof of Proposition 1.3 Lemma 3.3 also holds for $G_{\omega}$, where $\omega \in E(G)$. Our next result shows that, when $v=\lambda(\omega)$, the structure of $G_v$ can be determined completely.
\\ \\
\noindent {\bf Proposition 3.4.} \it  Let $\omega\in\Omega$ be an elliptic element,
and let let $v=\lambda(\omega)$ be its image under the building map.
There are two possibilities.
\begin{itemize}
\item[(i)] If $\omega\not\equiv\overline{\omega}$, then
$$G_v=G_{\omega} \cong \FF_{q^2}^*,$$ in which case $|\;G_v|=q^2-1$.
\item[(ii)] If $\omega \equiv \overline{\omega}$, then
$$G_v \cong GL_2(\FF_q),$$ in which case $|\;G_v|=q(q-1)^2(q+1)$.
\end{itemize}
\rm

\noindent {\bf Proof.} \rm
By Proposition 1.3 and Lemma 3.1, $G_{\omega}\cong \FF_{q^2}^*$
and $G_{\omega} \leq G_v$. Thus $G_{\omega}$ is completely reducible,
since $|G_{\omega}|$ is prime to $p$. Let $P$ be a matrix for which
$\widehat{G_v}=P^{-1}G_vP$ contains
$$\widehat{Z}=\{\diag(\alpha,\alpha^q):\alpha \in \FF_{q^2}^*\}.$$
\noindent It is clear that Lemmas 3.2 and 3.3 apply to the matrices
in $\widehat{G_v}$. Let
$$ M=\left[\begin{array}{lll} a & b\\[10pt]
c & d\end{array}\right]\in \widehat{G_v}.$$ \noindent By considering
the trace of the product $DM$, where $D=\diag(\alpha,\alpha^q)$, it
follows that
$$\alpha a+\alpha^q d \in \FF_q,\;\mathrm{for\; all}\; \alpha \in
\FF_{q^2}.$$ \noindent Now $\alpha(a+d) \in \FF_{q^2}$. From the
case where $\alpha \neq \alpha^q$ we deduce that $a,d \in
\FF_{q^2}$. \\ \\
\noindent Now $\alpha a+ \alpha^q a^q \in \FF_q$ and so
$$\alpha^q(d-a^q) \in \FF_q,\; \mathrm{for\; all}\; \alpha \in
\FF_{q^2}.$$ \noindent We conclude that $d=a^q$ and hence that $bc
\in \FF_q$. \noindent Suppose that $bc \in \FF_q^*$. Then by Lemma
3.3
$$ N=\left[\begin{array}{lll} 0 & b\\[10pt]
c & 0\end{array}\right]\in \widehat{G_v}.$$ \noindent Let
$$ N_1=\left[\begin{array}{lll} 0 & b_1\\[10pt]
c_1 & 0\end{array}\right]\in \widehat{G_v}.$$ \noindent  Now $NN_1
\in \widehat{G_v}$. It follows from the above that $bc_1,b_1c \in
\FF_{q^2}$ with $(bc_1)^q=b_1c$ and $bc,b_1c_1 \in \FF_q^*$. We
deduce that
$$b_1=b\beta \;\mathrm{and}\;c_1=c\beta^q ,$$
for some $\beta \in \FF_{q^2}$, i.e. $N_1=\diag(\beta,\beta^q)N$.
\\ \\
\noindent There are two possibilities. If $G_v=G_{\omega}$, then
$G_v\cong \FF_{q^2}^*$ and we are finished. Suppose then that there
exists
$$\left[\begin{array}{lll} a & b\\[10pt]
c & d\end{array}\right]\in \widehat{G_v},$$ with $b \neq 0$ or
$c\neq 0$. We now show that in this case $bc \neq 0$. \noindent If
$b\neq 0$ and $c=0$ or vice versa, then by Lemma 3.3
$$ S=\left[\begin{array}{lll} 1 & b\\[10pt]
0 &1\end{array}\right]\in \widehat{G_v}\;\;\mathrm{or}\;\;T= \left[\begin{array}{lll} 1 & 0\\[10pt]
c & 1\end{array}\right]\in \widehat{G_v}.$$
\noindent Let $G_{\omega}=\langle M_0\rangle$ and
$P^{-1}M_0P=D=\diag(\gamma,\gamma^q)$, where $\gamma \neq \gamma^q$.
Now $M_0$ fixes (distinct) elliptic points $\omega,
\overline{\omega}$, say. Then $D$ fixes $0,\infty$ and we may assume
that $P(\infty)=\omega$ and $P(0)=\overline{\omega}$. Now
$S(\infty)=\infty$. If $S \in \widehat{G_v}$, then $PSP^{-1} \in
G_{\omega}$. But by Proposition 1.3 $G_{\omega}$ does not contain
any elements of order $p$. Hence $\widehat{G_v}$ contains no
matrices of type $S$ and similarly none of type $T$. \noindent We
may assume therefore that $bc=\varepsilon \in \FF_q^*$. Let $D,N$ be as
above. Now $\{1,\gamma\}$ is a basis of $\FF_{q^2}$ over $\FF_q$. It
follows from Lemma 3.3 and the above that $\widehat{G_v}$ is the
subset of invertible elements in the $4$-dimensional central simple
algebra over $\FF_q$ with $\FF_q$-basis $\{I_2,D,N,DN\}$. It is
clear that
$$|\widehat{G_v}|=|\{(\rho,\epsilon) \in \FF_{q^2} \times
\FF_{q^2}:\rho\rho^q\neq\varepsilon
\epsilon\epsilon^q\}|=q(q-1)^2(q+1).$$ \noindent By Wedderburn's
Theorem it follows that
$$\widehat{G_v} \cong GL_2(\FF_q).$$
\noindent Now $N(\infty)=0$ and $N(0)=\infty$ and so $PNP^{-1}\in G_v$
``interchanges" $\omega$ and $\overline{\omega}$.
(Note that if $g(\omega)=\overline{\omega}$ then $g(\overline{\omega})=\omega.)$
\hfill $\Box$
\\ \\
Our next lemma highlights the importance of the $q^2-1$ as a feature of our results.
\\ \\
\noindent {\bf Lemma 3.5.} \it For a vertex stabilizer $G_v$ the following three statements
are equivalent:
\begin{itemize}
\item[(i)] $q^2 -1$ divides $|G_v|$.
\item[(ii)] $G_v$ contains a matrix whose eigenvalues are not in $\FF_q$.
\item[(iii)] $G_v$ contains a cyclic subgroup of order $q^2-1$.
\end{itemize}
\rm \noindent
{\bf Proof.}
$(i)\Rightarrow (ii)$:
If $q+1$ is divisible by an odd prime $r$, then $r$ divides neither
$q$ nor $q-1$. Let $M\in G_v$ be an element of order $r$. From
the order we see that the eigenvalues of $M$ cannot be in $\FF_q$.
If $q+1$ is not divisible by any odd prime, then it is divisible by $4$.
If $G_v$ contains an element of order $4$, we can argue as before.
If not, we fix a $2$-Sylow subgroup of $P$ of $G_v$, which then is
necessarily of exponent $2$ and hence abelian. So all matrices in
$P$ can be simultaneously diagonalized. Since the eigenvalues can
only be $1$ and $-1$, there are only $4$ such diagonal matrices.
But the order of $P$ is divisible by $8$, a contradiction.
\par
$(ii)\Rightarrow (iii)$:
If the eigenvalues of $M$ are not in $\FF_q$, then by Lemma 3.2 they
are in $\FF_{q^2}\setminus\FF_q$. Let
$$I(M)= \left\{\alpha I_2+\beta M: \alpha,\beta \in \FF_q^*,\;(\alpha,\beta) \neq (0,0)\right\}.$$
\noindent By Lemma 3.3 then $I(M) \leq G_v$. Part (iii) follows since $I(M) \cong \FF_{q^2}^*$.
\par
$(iii)\Rightarrow (i)$ is trivial.
\hfill $\Box$
\\ \\
\noindent For some $q$ (for example $q=4$) every subgroup of $G$ of order $q^2-1$ is cyclic. On the other hand for the case $q=3$ the embedding of $A_4$ in $PGL_2(\FF_3)$ gives rise to a subgroup $S$ of $G$ containing $Z$ of order $8$ for which $S/Z$ is not cyclic.
\par
In Lemma 3.5 the condition (i) can be replaced by
$$(i)'\;\;|G_v|\;is\;divisible\;by\;q+1\;(q\neq 3)\;and\;8\;(q=3).$$
\noindent Here the restriction when $q=3$ is necessary. It is well-known [Se, p.86] that, when $A=\FF_3[t]$,
there is a vertex $v'$ for which
\begin{itemize}
\item[(1)]$|G_{v'}|=12$,
\item[(2)] every matrix in $G_{v'}$ has eigenvalues in $\FF_3^*$.
\end{itemize}
\noindent We note that the proof of Lemma 3.5 shows that when, $q \neq 3$, the following implication holds.
$$  q+1\;divides\;|G_v| \Rightarrow q^2-1\; divides \; |G_v|.$$
\noindent
As stated above the main aim in this section is to prove that the converse
of Proposition 3.4 holds. We will prove that, if $q^2-1$ divides $|G_v|$,
then $v=\lambda(\omega)$ for some elliptic point $\omega \in E(G)$.
We require one more lemma.
\\ \\
\noindent {\bf Lemma 3.6.} \it Let $\delta$ be odd and let $M\in G$ be
a matrix of finite order whose eigenvalues are not in $\FF_q$. Then
\begin{itemize}
\item[(i)]  $M$ does not fix any edges of $\TTT$.
\item[(ii)] $M$ fixes exactly one vertex of $\TTT$.
\end{itemize}
\rm

\noindent
{\bf Proof.}
(i) Suppose that $M$  fixes an edge. Then there exists a matrix
$P\in GL_2(K_{\infty})$ that maps this edge to the standard edge
whose stabilizer is $Z_\infty\cdot {\cal J}$ where $Z_\infty$ is the centre
of $GL_2(K_\infty)$ and ${\cal J}$ is the Iwahori group
$${\cal J}=\left\{{a\ b\choose c\ d}\in GL_2(\OOO_\infty)\ :\ c\in\pi\OOO_\infty\right\}.$$
From the determinant we see that $P$ conjugates $M$ into ${\cal J}$.
Let $\widetilde{M}\in {\cal J}$ be this conjugate of $M$. Then the
characteristic polynomial $X^2 -\tau X+\eta$ of $\widetilde{M}$ is
irreducible over $\FF_q$, and hence over
$\FF_{q^{\delta}}$ if $\delta$ is odd.
\par
On the other hand, reducing $\widetilde{M}$ modulo the maximal
ideal of $\OOO_{\infty}$ we obtain a matrix with the same characteristic
polynomial. But the reduced matrix has the form
${a\ b\choose 0\ d}$ with entries in $\FF_{q^{\delta}}$.
So its characteristic polynomial splits over $\FF_{q^{\delta}}$,
a contradiction.
\par
(ii) By [Se, Proposition 2, p.79] $M$ fixes at least one vertex.
If $M$ fixes two different vertices of $\TTT$, then it fixes the whole
geodesic on $\TTT$ between these two vertices and hence at least
one edge in contradiction to (i).
\hfill $\Box$
\\ \\
By the way, Lemma 3.6 (ii) provides an alternative proof of the claim in
Lemma 3.1 that $\lambda(\omega)=\lambda(\overline{\omega})$.
\\ \\
\noindent We now come to the principal results of this section.
\\ \\
\noindent {\bf Theorem 3.7.} \it Let $\delta$ be odd and let $ v \in \vert(\TTT)$. Then
$$ v = \lambda(\omega),\; for\; some\; \omega \in E(G),\;\; if\;\; and\;\; only\;\; if\;\; q^2-1\; divides\; |G_v|.$$

\noindent {\bf Proof.} \rm If $v=\lambda(\omega)$ for some $\omega\in E(G)$, then $q^2 -1$ divides
$|G_v|$ by Proposition 3.4.
\par
Conversely, assume that $q^2 -1$ divides $|G_v|$. Then $G_v$ contains a cyclic subgroup $C$ of order
$q^2 -1$ by Lemma 3.5. By Lemma 1.6 this subgroup $C$ fixes an elliptic point $\omega\in E(G)$.
Again by Proposition 3.4 we know that $q^2 -1$ divides $|G_{v'}|$ for $v'=\lambda(\omega)$. Since
$C$ is contained in $G_v$ and $G_{v'}$, Lemma 3.6 implies $v'=v$.
\hfill $\Box$
\\ \\
\noindent {\bf Theorem 3.8.} \it If $\delta$ is odd, there exist natural bijections between the
following sets
\begin{itemize}
\item[(i)] vertices $\widetilde{v}$ of $G\!\setminus\!\TTT$ such that $q^2 -1$ divides $|G_v|$;
\item[(ii)] conjugacy classes (in $G$) of cyclic subgroups of $G$ of order $q^2 -1$;
\item[(iii)] the orbits of the $Gal(\widetilde{K}/K)$-action on $\Ell(G)$.
\end{itemize}
\rm

\noindent {\bf Proof.} \rm We first establish the bijection between (i) and (ii).
Let $v$ be a vertex of $\TTT$ with image $\widetilde{v}$ in $G\!\setminus\!\TTT$.
If $G_v\cong\FF_{q^2}^*$, this is such a cyclic subgroup of order $q^2 -1$, and
the stabilizers of the other lifts of $\widetilde{v}$ to $\vert(\TTT)$ are exactly
the conjugates of $G_v$. A similar argument applies if $G_v \cong GL_2(\FF_q)$.
Of course, then $G_v$ has several cyclic subgroups of order $q^2 -1$, but they are all
conjugate (already in $G_v$).
\par
Conversely, let $C$ be a cyclic subgroup of $G$ of order $q^2 -1$. By Lemma 3.6 it fixes
exactly one vertex of $\TTT$. So its conjugacy class fixes exactly one vertex of
$G\!\setminus\!\TTT$.
\par
The bijection between (ii) and (iii) follows by applying the action of $G$ to the bijection
in Lemma 1.6.
\hfill $\Box$
\\ \\
\noindent {\bf Remark 3.9.}
Theorem 3.8 (in combination with Proposition 3.4) implies in particular that over
every vertex $\widetilde{v}$ of $G\!\setminus\!\TTT$ with $G_v\cong GL_2(\FF_q)$
there lies exactly one elliptic point of $G\!\setminus\!\Omega$; and over every
vertex $\widetilde{v}$ of $G\!\setminus\!\TTT$ with $G_v\cong\FF_{q^2}^*$ lie two ($Gal(\widetilde{K}/K)$-conjugate) elliptic points of $G\!\setminus\!\Omega$.
\par
But when considering the building map $\lambda:\Omega\to\TTT$, over every
vertex $v$ of $\TTT$ with $G_v\cong GL_2(\FF_q)$ there lie $q(q-1)$ elliptic
points on $\Omega$, in $q(q-1)/2$ pairs of $Gal(\widetilde{K}/K)$-conjugate
elliptic points, corresponding to the $q(q-1)/2$ different cyclic subgroups
of order $q^2 -1$ in $GL_2(\FF_q)$. (Compare Lemmas 1.6 and 3.6.) Over every
vertex $v$ of $\TTT$ with $G_v\cong\FF_{q^2}^*$ we again have one pair of
$Gal(\widetilde{K}/K)$-conjugate elliptic points on $\Omega$.
\par
One should not forget however that there also are uncountably many
non-elliptic points lying over each of these vertices, as for every
vertex $v$ of $\TTT$ there are uncountably many points of $\Omega$
mapping to $v$ under the building map.
\\ \\
A much more general statement than Proposition 3.4, namely the complete
classification of all possible types of vertex stabilizers for any constant
field (not just for $\FF_q$) and for any $\delta$ is given in [MS3].
\\

\subsection*{4. Isolated vertices and amalgams}

A vertex $\widetilde{v}$ of the quotient graph $G\backslash\TTT$ is called {\it isolated}
if there is only one edge of $G\backslash\TTT$ attached to it. Obviously this is equivalent
to $G_v$ acting transitively on the $q^{\delta}+1$ edges of $\TTT$ attached to $v$.
\\ \\
\noindent {\bf Theorem 4.1.} \it Let $ v\in \vert(\TTT)$. Then
$\widetilde{v}$  is an isolated vertex of $G\backslash\TTT$ if and
only if the following two conditions both hold:
\begin{itemize}
\item[(i)] $\delta=1$,
\item[(ii)] $G_v$ satisfies any of the three equivalent conditions of Lemma 3.5.
\end{itemize}

\noindent {\bf Proof.} \rm
Assume first that $\delta=1$ and $G_v$ contains a cyclic group of
order $q^2 -1$. Then by Lemma 3.6 none of the elements outside
$Z$ can fix an edge. So $G_v$ acts transitively on the $q+1$ edges
adjacent to $v$, and $\widetilde{v}$ is isolated.
\par
Now assume conversely that $\widetilde{v}$ is isolated. Then $G_v$
acts transitively on the $q^{\delta}+1$ edges emanating from $v$.
So $|G_v|$ is divisible by $(q-1)(q^{\delta}+1)$.
\par
If $q^{\delta}+1$ is divisible by an odd prime $r$, then $r$ divides
neither $q$ nor $q-1$. Let $M\in G_v$ be an element of order $r$.
From the order we see that the eigenvalues of $M$ cannot be in
$\FF_q$. But by Lemma 3.2 they are in $\FF_{q^2}$, so $r$ divides
$q+1$. Together with $r$ dividing $q^{\delta}+1$ this implies that
$\delta$ is odd. If $\delta$ were bigger than $1$, then $G_v$
would act transitively on at least $q^3 +1$ edges. But
$|G_v/Z|\leq q^3 -q$ by Proposition 3.4.
\par
If $q^{\delta}+1$ is not divisible by any odd prime, then it is divisible
by $4$, and hence $q$ is congruent to $3$ modulo $4$ and $\delta$ is
odd. As above we obtain $\delta=1$. Moreover, $q+1$ divides $|G_v|$
because it divides $q^{\delta}+1$.
\hfill $\Box$
\\ \\
Combining Theorem 4.1 with Theorem 3.8 we obtain the following
\\ \\
\noindent {\bf Corollary 4.2.} \it
Let $\delta=1$. Then the building map induces a bijection between the
$Gal(\widetilde{K}/K)$-orbits on elliptic points of $G\setminus\Omega$
and the isolated vertices of $G\setminus\TTT$ with the properties
described in Proposition 3.4.
\par
The number of isolated vertices with stabilizer isomorphic to
$GL_2(\FF_q)$ (resp. to $\FF_{q^2}^*$) is $|\Ell(G)_2|$ (resp.
$r=\frac{1}{2}(|\Ell(G)|-|\Ell(G)_2|)$. In particular, these
numbers only depend on $K$, not on the choice of the degree
one place $\infty$.
\rm
\\ \\
\noindent We also record a graph-theoretic property of isolated
vertices.
\\ \\
\noindent {\bf Proposition 4.3.} \it
\begin{itemize}
\item[a)] Let $\delta$ be odd and let $v_1,v_2 \in \vert(\TTT)$,
 where $|G_{v_i}|$ is divisible by $q^2-1$, $(i=1,2)$.
 Then the (geodesic) distance between $v_1$ and $v_2$ (in $\TTT$)
and consequently the distance between $\widetilde{v_1}$ and
$\widetilde{v_2}$ (in $G\backslash \TTT$) is even.
\item[b)] The distance between any two isolated vertices of
$G\!\setminus\!\TTT$ is even.
\end{itemize}

\noindent {\bf Proof.} \rm a) By Theorem 3.7 there exist $\omega_i
\in E(G)$ with $v_i =\lambda(\omega_i)$, ($i=1,2$). Fix
$\varepsilon\in\FF_{q^2}\setminus\FF_q$. By Theorem 1.4 we can write
$\omega_i =\frac{\varepsilon +s_i}{t_i}$ with $s_i ,t_i \in A$. Thus
$\omega_2 =M(\omega_1)$ with
$$M=\left[\begin{array}{cc} 1 & s_2\\[10pt]
0 & t_2 \end{array}\right]\left[\begin{array}{cc} t_1 & -s_1\\[10pt]
0 & 1 \end{array}\right]\in GL_2(K_\infty).$$ Consequently $v_2
=M(v_1)$ by [GR, (1.5.3)]. Let $d(v_1,v_2)$ be the distance between
$v_1$ and $v_2$. Then by [Se, Corollary, p.75] and Lemma 1.7 b)
$$d(v_1,v_2) \equiv \nu(\det(M))\equiv \nu(t_1)+\nu(t_2)\equiv 0\;(\mod 2).$$
\noindent Part (b) follows from part (a) and Theorem 4.1.
\hfill $\Box$
\\ \\
\noindent The principal group-theoretic consequence of Theorem 4.1
is the following.
\\ \\
\noindent {\bf Theorem 4.4.} \it Suppose that $\delta=1$ and that
$|G_v|$ is divisible by $q^2-1$. There are two possibilities.
\begin{itemize}
\item[(i)] If $G_v \cong GL_2(\FF_q)$, then there exists a subgroup
$H$ of $G$ such that
$$G \cong GL_2(\FF_q)\star_{\quad B_2(\FF_q)}H,$$
\noindent where $B_2(\FF_q)$ is the usual Borel subgroup (of order
$q(q-1)^2$).
\item[(ii)] If $G_v\cong \FF_{q^2}^*$, then there exists a subgroup
$H$ of $G$ for which
$$G \cong \FF_{q^2}^*\star_{\quad Z }H.$$
\noindent Hence
$$PGL_2(A) \cong (\ZZ/(q+1)\ZZ)\star H',$$
\noindent where $H'=H/Z$.
\end{itemize}

\noindent {\bf Proof.} \rm Let $e$ be any edge incident with $v$. Then
by Theorem 4.1 $\widetilde{v}$ is isolated (in $G \backslash\TTT$) and
so $|G_v:G_e|=q+1$. Bass-Serre theory [Se, Theorem 13, p.55]
presents $G$ as the {\it fundamental group of a graph of groups}
[Se, p.42] given by a {\it lift} $$j: \TTT_0\longrightarrow \TTT,$$
\noindent where $\TTT_0$ is a maximal subtree of $G\backslash \TTT$.
Since $\widetilde{v}$ is isolated $v,e$ may be assumed to lie in $j(\TTT_0)$.
It follows that
$$G \cong G_v \star_{\quad G_e}H,$$
where $H \neq G_v,G_e$ is some subgroup of $G$. The results follow.
\hfill $\Box$
\\ \\
\noindent When $\delta=1$ a decomposition of type (i) always occurs
because the standard vertex has stabilizer $GL_2(\FF_q)$. More
interesting decompositions occur when there are isolated vertices of
type (ii).
\\ \\
\noindent {\bf Theorem 4.5.} \it Suppose that $\delta=1$. Then there
exists a subgroup $P$ of $PGL_2(A)$ for which the following free
product decomposition holds
$$PGL_2(A)\cong\left(\Star_{i=1}^r\ZZ/(q+1)\ZZ\right)\star P,$$
\noindent where
$$2r=|\Ell(G)|-|\Ell(G)_2|.$$
Moreover $r$ is maximal in the following sense. Suppose that $C$ is a cyclic subgroup of $PGL_2(A)$
of order $q+1$ for which
$$PGL_2(A)=C \star Q.$$
\noindent Then there exists $v \in \vert(\TTT)$ such that
\begin{itemize}
\item[(i)] $G_v \cong \FF_{q^2}$,
\item[(ii)] $\psi(G_v)=C$,
where $\psi: G \rightarrow  PGL_2(A)$ is the natural map.
\end{itemize}

\noindent {\bf Proof.} \rm Let
$$\widetilde{V}=\{ \widetilde{v} \in \vert(G \backslash \TTT):\; G_v
\cong \FF_{q^2}^* \}.$$ \noindent Let
$\widetilde{v_1},\;\widetilde{v_2} \in \widetilde{V}$. Then, by
Theorem 3.7, $G_{v_i}=G_{\omega_i}$, for some $\omega_i \in E(G)$,
where $\omega_i \not\equiv \overline{\omega_i}$, ($i=1,2$). If
$\widetilde{v_1}=\widetilde{v_2}$, then $v_2=g(v_1)$, for some $g
\in G$, so that $G_{\omega_2}=gG_{\omega_1}g^{-1}=G_{g(\omega_1)}$.
It follows that
$\{\omega_2,\;\overline{\omega_2}\}=\{g(\omega_1),\;g(\overline{\omega_1})\}$.
\noindent On the other hand if $\omega_j (\neq \overline{\omega_j})
\in E(G)$ and $S_j=\{\omega_j\;\overline{\omega_j}\}$, where
$j=3,4$, then for all $g \in G$ either $S_3=g(S_4)$ or $S_3 \cap
g(S_4) = \emptyset$. By Corollary 4.2 we have
$|\widetilde{V}| =r$, where $r$ is defined as above.
\par
The free product decomposition is a consequence of an iteration
of the process described in the proof of Theorem 4.4 (ii).
\par
For the last part of the theorem $C$, under $\psi$, lifts to a cyclic subgroup $C'$
of $G$ of order $q^2-1$. Now by [Se, Proposition 2, p.76] $C' \leq G_v$, for some $v \in \vert( \TTT)$.
Then by Theorem 3.7 there are two possibilities for $G_v$, described in
Proposition 3.4. Either $G_v=C'$ in which case we are finished, or $G_v\cong PGL_2(\FF_q)$. In the latter case
the canonical map from $PGL_2(A)$ onto $C$ restricts to an epimorphism
$$PGL_2(\FF_q)\twoheadrightarrow C.$$
\noindent This gives the desired contradiction.
\hfill $\Box$
\\ \\
\noindent {\bf Theorem 4.6.} \it
\begin{itemize}
\item[(a)] For $q\geq 8$ and $g>0$ there exist exactly two rings
$A$ (up to isomorphism) such that all isolated vertices of
$GL_2(A)\!\setminus\!\TTT$ have stabilizers isomorphic to
$GL_2(\FF_q)$, namely $A=\FF_9[x,y]$ with $y^3+y=x^4$ (genus $3$)
or with $y^2=x^3-x$ (genus $1$).
\item[(b)] For fixed $q\geq 8$ the number $r$ of free factors in
Theorem 4.5 grows exponentially with $g$.
More precisely, for $q\geq 8$ and all cases of positive genus
except the two discussed in part (a) we have
$$r\geq\frac{1}{4}(\sqrt{q}-1)^{2g}>\frac{3^g}{4}.$$
\end{itemize}
\rm

\noindent {\bf Proof.}
(a) From Theorem 2.9 (a) we know already that there are only two
fields $K$ with these properties. For any choice of the place
$\infty$ of degree $1$ we get a ring $A$ with this property.
It remains to show that different choices of $\infty$ give
isomorphic rings.
\par
For the genus $3$ case we use that
by [St, Exercise 6.10] the automorphism group of a Hermitian
function field acts transitively on its places of degree $1$. So
different choices of $\infty$ will lead to isomorphic rings $A$.
\par
The elliptic case can be seen by some easy calculations with
Weierstrass equations.
\\
(b)
If $\Ell(G)_2$ is strictly smaller than $\Ell(G)$, then, because of
the group structure, it has index at least $2$. Hence, if there are
elements of order bigger than $2$ in $\Ell(G)$, their number is at
least $\frac{1}{2}L_K(-1)\geq\frac{1}{2}(\sqrt{q}-1)^{2g}$. So in
that case the number of isolated vertices with cyclic stabilizer
is at least $\frac{1}{4}(\sqrt{q}-1)^{2g}$, which for $q\geq 8$ is
bigger than $\frac{1}{4}3^g$.
\hfill $\Box$
\\ \\
\noindent {\bf Remarks 4.7.}
\begin{itemize}
\item[(a)] It is, of course, well possible that the group $PGL_2(A)$ also
splits off other free factors than those stipulated by Theorem 4.5. Let
for example $A=\FF_9[x,y]$ with $y^2 =x^3 -x$. Then $r=0$, but from
Takahashi's results [T] one obtains that in this case $PGL_2(A)$ is
a free product of $10$ infinite groups.
\item[(b)] By the same arguments as in the proof of Theorem 4.6 (b),
for $q\in\{5,7\}$ we still have $r>\frac{1}{4}(\frac{3}{2})^g$
{\it provided} $r$ is not zero. (Compare Remark 2.10.)
\item[(c)] Theorem 4.5 has a number of interesting consequences.
For example suppose that $r\geq 2$ and that $q \equiv -1\;(\mod 6)$.
Then there exists an epimorphism
$$\theta: PGL_2(A) \twoheadrightarrow PSL_2(\ZZ),$$
\noindent since
$$PSL_2(\ZZ) \cong (\ZZ/2\ZZ)\star (\ZZ/3\ZZ).$$
\end{itemize}

\noindent {\bf Example 4.8.}
Let
$$A=\FF_2[x,y]\ \ \ \hbox{\rm with}\ \ \  y^2 +y=x^3 +x+1.$$
This elliptic curve has exactly one rational point, namely the
one at infinity. So $L_K(u)=1-2u+2u^2$, and thus $L_K(-1)=5$ and $r=2$.
More precisely,
$$PGL_2(A)\cong GL_2(A)\cong\ZZ/3\ZZ\star\ZZ/3\ZZ\star\Delta(\infty),$$
where $\Delta(\infty)=B_2(A)\star_{B_2(\FF_2)} GL_2(\FF_2)$
(cf. Takahashi [T] and [MS1, Theorem 5.3] or (the proof of)
[MS2, Lemma 5.2 (c)]).
Since the normal hull of $B_2(A)$ in $GL_2(A)$ contains all
elements from $\Delta(\infty)$, we see that a finite group can
be generated by two elements of order $3$ if and only if it is
the quotient of this $GL_2(A)$ by a normal non-congruence
subgroup of level $A$.
For results on which classical finite simple groups can be
generated by two elements of order $3$ see [LS, Corollary 1.8].
\\ \\
{\bf Example 4.9.}
Let
$$A=\FF_7[x,y]\ \ \ \hbox{\rm with}\ \ \ y^2 =x^3 +4.$$
Then $L_K(u)=1-5u+7u^2$. Thus $L_K(-1)=13$ and $r=6$.
So there exists a surjective homomorphism from $GL_2(A)$
to any finite (or infinite) group that is generated by
at most $6$ elements of orders dividing $8$.
More precisely, by Takahashi's description of the quotient
graph (cf. [T]) we have
$$PGL_2(A)\cong\left(\Star_{i=1}^6\ZZ/8\ZZ\right)\star
\Delta(0)\star\Delta(\infty),$$
where $\Delta(0)$, $\Delta(\infty)$ are infinite subgroups
and, again, $\Delta(\infty)$ contains all upper triangular
matrices (modulo $Z$).
\par
In particular, there exists a normal non-congruence subgroup
$N$ of level $A$ such that $G/N$ is isomorphic to the
permutation group of Rubik's cube (which is generated by
$6$ elements of order $4$). Recall that the order of
that permutation group is roughly $43\cdot 10^{18}$.
\\ \\ \\
{\bf Acknowledgements.} This paper is part of a project supported by grant
99-2115-M-001-011-MY2 from the National Science Council (NSC) of Taiwan.
The biggest part was written while the second author was working at the
Institute of Mathematics at Academia Sinica in Taipei. He wants to thank
Julie Tzu-Yueh Wang, Wen-Ching Winnie Li, Jing Yu, Liang-Chung Hsia and
Chieh-Yu Chang for help in general as well as for help with the application
for that grant. During the final stage the second author was supported by
ASARC in South Korea.
\par
Several ideas in the paper were developed during research visits of the
second author at Glasgow University. The hospitality of their Mathematics
Department is gratefully acknowledged.
Finally, the second author thanks Ernst-Ulrich Gekeler for helpful
discussions during a research visit to Saarbr\"ucken.
\\

\subsection*{\hspace*{10.5em} References}
\begin{itemize}
\item[{[BD]}] V.~Bautista-Ancona and J.~Diaz-Vargas: Quadratic
function fields with exponent two ideal class group, \it
J. Number Theory \bf 116 \rm (2006), 21-41
\item[{[B]}] N.~Bourbaki: \it Commutative Algebra,
\rm Addison-Wesley, London, 1972
\item[{[D]}] V.~G.~Drinfeld: Elliptic Modules, \it Math. USSR-Sbornik
\bf 23 \rm (1976), 561-592
\item[{[G]}] E.-U.~Gekeler: \it Drinfeld Modular Curves, \rm
Springer LNM 1231, Berlin Heidelberg New York, 1986
\item[{[GR]}] E.-U.~Gekeler and M.~Reversat: \rm Jacobians of Drinfeld
Modular Curves, \it J. Reine Angew. Math. \bf 476 \rm (1996), 27-93
\item[{[LS]}] M.~Liebeck and A.~Shalev: Classical groups, probabilistic
methods, and the $(2,3)$-generation problem, \it Ann. of Math. \bf 144
\rm (1996), 77-125
\item[{[M]}] A.~W.~Mason: Serre's generalization of Nagao's
theorem: an elementary approach, \it Trans. Amer. Math. Soc.
\bf 353 \rm (2003), 749-767
\item[{[MS1]}] A.~W.~Mason and A.~Schweizer: The minimum index of a
non-congruence subgroup of $SL_2$ over an arithmetic domain. II: The
rank zero cases, \it J. London Math. Soc. \bf 71 \rm (2005), 53-68
\item[{[MS2]}] A.~W.~Mason and A.~Schweizer: Non-standard automorphisms
and non-cong\-ruence subgroups of $SL_2$ over Dedekind domains contained
in function fields, \it J. Pure Appl. Algebra \bf 205 \rm (2006), 189-209
\item[{[MS3]}] A.~W.~Mason and A.~Schweizer: The stabilizers in a Drinfeld
modular group of the vertices of its Bruhat-Tits tree: an elementary approach,
\it submitted, \rm http://arxiv.org/abs/1203.3617
\item[{[R1]}] M.~Rosen: The Hilbert class field in function fields,
\it Expo. Math. \bf 5 \rm (1987), 365-378
\item[{[R2]}] M.~Rosen: \it Number Theory in Function Fields,
\rm Springer GTM 210, Berlin Heidelberg New York, 2002
\item[{[RSt]}] H.-G.~R\"uck and H.~Stichtenoth: A characterization
of Hermitian function fields over finite fields,
\it J. Reine Angew. Math. \bf 457 \rm (1994), 185-188
\item[{[Se]}] J.-P.~Serre: \it Trees, \rm Springer, Berlin
Heidelberg New York, 1980
\item[{[St]}] H.~Stichtenoth: \it Algebraic Function Fields and
Codes (Second Edition), \rm Springer GTM 254, Berlin Heidelberg,
2009
\item[{[T]}] S.~Takahashi: The fundamental domain of the
tree of $GL(2)$ over the function field of an elliptic curve,
\it Duke Math. J. \bf 72 \rm (1993), 85-97
\item[{[ZS]}] O.~Zariski and P.~Samuel: \it Commutative Algebra
(Volume 1), \rm Springer GTM 21, Berlin Heidelberg New York, 1975
\end{itemize}

\end{document}